\documentclass[A4paper,12pt]{revtex4}

\usepackage{epsfig}
\usepackage{subfigure}
\begin{document}

\title{Smoothing Noisy Spectroscopic Data with Many-Knot Spline Method}

\author{M.H.Zhu}
\thanks{peter\_zu@163.com.}
\affiliation{Space Exploration Laboratory, Macao University of
Science and Technology, Taipa, Macao}
\author{L.G.Liu, D.X.Qi, Z.You, A.A.Xu}
\affiliation{Faulty of Information Technology, Macao University of
Science and Technology, Taipa, Macao}

\begin{abstract}
In this paper, we present the development of a many-knot spline
method derived to remove the statistical noise in the
spectroscopic data. This method is an expansion of the B-spline
method. Compared to the B-spline method, the many-knot spline
method is significantly faster.
\end{abstract}


\maketitle

Smoothing using least square method with B-spline functions is
helpful in reducing the statistical noise in the spectroscopic
data, such as gamma-ray spectrum. However, when the amount of
points is very large, this method becomes time-consuming because
of the need to solve nonlinear least squares equations, especially
if the initial knots are not well determined or too many channels
are selected as initial knots.
\begin{figure}[h]
  \centering
  \subfigure[]{
    \label{fig:subfig:a}
    \includegraphics[height=6cm,width=7cm]{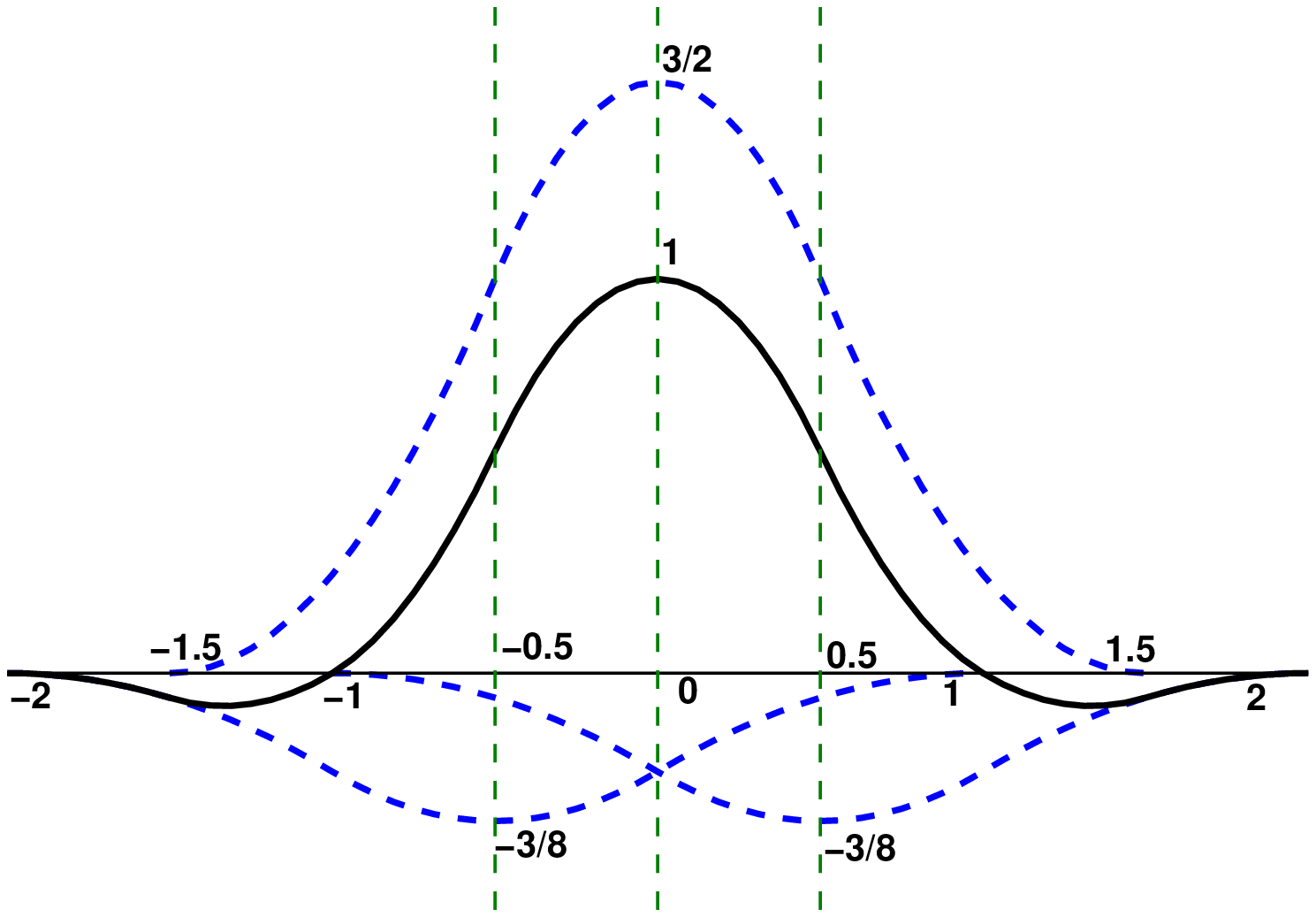}}
    \hspace{.2in}
  \subfigure[]{
    \label{fig:subfig:b}
    \includegraphics[height=6cm,width=7cm]{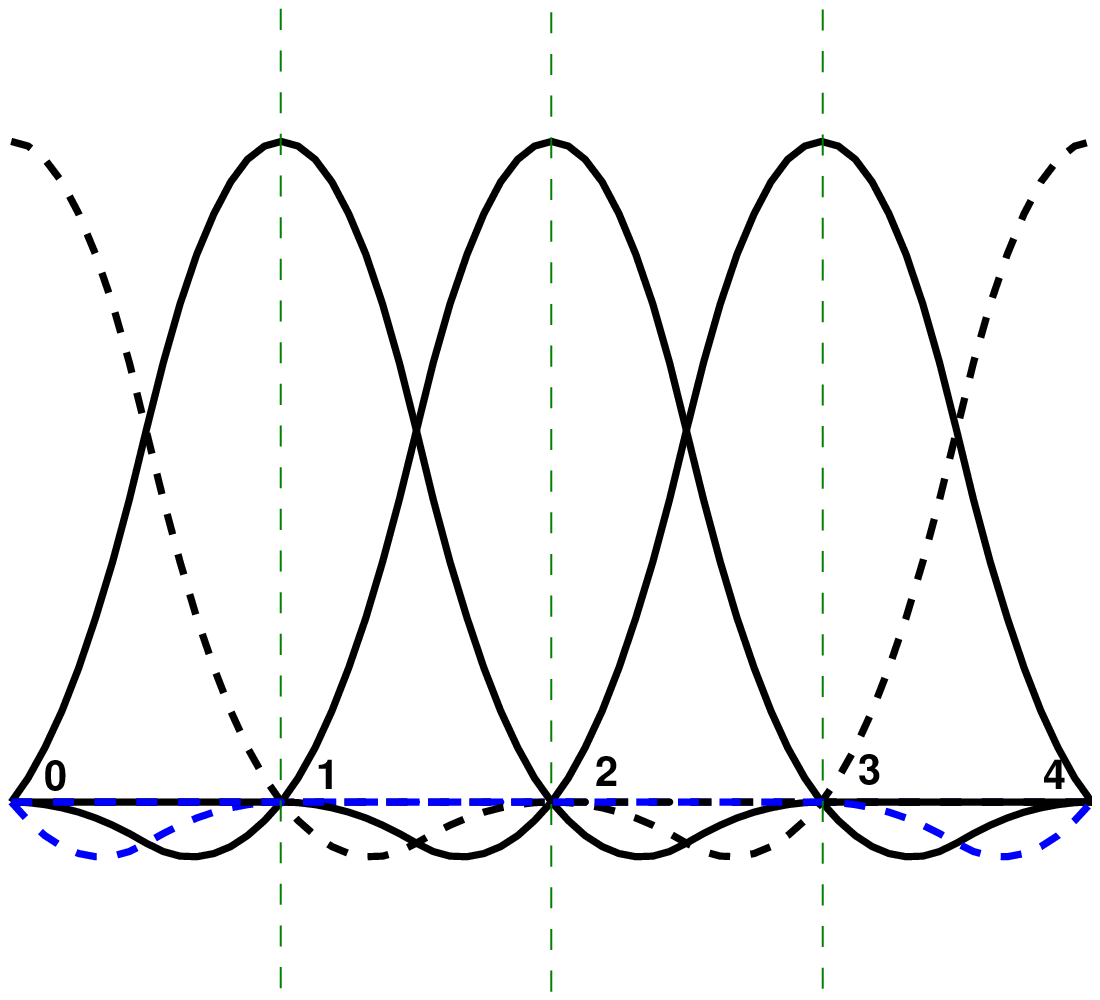}}
  \caption{(a) The construction of quadric many-knot spline basis
function shape; (b) 5-basis many-knot splines with each count
value = 1.}
  \label{figure1}
\end{figure}

To avoid this problem, a many-knot spline is used in substitution
for B-spline for smoothing noisy spectroscopic data. The many-knot
method is an expansion technique of B-spline, described by the
third author D.X.Qi\cite{1,2}. Its basis function of order
\emph{\textbf{k}} is denoted as following:
\def\formula1{
   q_k (x) = \sum\limits_{i = 0}^{k - 1} {t_i } \Omega _k^{ < a_i  >
}    }
  \overfullrule 5pt
  \mathindent\linewidth\relax
  \advance\mathindent-284pt
  \begin{equation}
    \label{eq:1}
    \formula1
  \end{equation}
where $t_{i}$ are coefficients that can be solved easily according
to its characteristic and $\Omega_{k}^{<l>}$ is
\def\formula1{
   \Omega _k^{ < l > }  = {\textstyle{1 \over 2}}[\Omega _k^{ <
x + l
> }  + \Omega _k^{ < x - l > } ],\quad l \ne 0;
    }
  \overfullrule 5pt
  \mathindent\linewidth\relax
  \advance\mathindent-284pt
  \begin{equation}
    \label{eq:1}
    \formula1
  \end{equation}
$a_{0}$,$a_{1}$,$...$,$a_{k-1}$ are transform coefficients and
$\Omega_{k}(x)$ is B-spline basis function of order
\emph{\textbf{k}}. Using $\emph{\textbf{k}}=2$ as an example, the
quadric many-knot spline basis function can be expressed as:
\def\formula1{
   q_2 (x) = 2\Omega _2 (x) - {\textstyle{1 \over 2}}[\Omega _2 (x +
{\textstyle{1 \over 2}}) + \Omega _2 (x - {\textstyle{1 \over
2}})]
    }
  \overfullrule 5pt
  \mathindent\linewidth\relax
  \advance\mathindent-284pt
  \begin{equation}
    \label{eq:1}
    \formula1
  \end{equation}
with its structure shown as Fig.1(a). Many-knot spline basis
functions have the similar characteristics as B-spline basis
functions. Importantly, unlike the B-spline method, the many-knot
spline method can obtain the best fitting curve without solving
the equation, which can save a lot of time in the calculation. The
details can be obtained from Ref.\cite{1,2,3}.

In the new method, five channels with equal interval in the
interesting region are selected as initial knots. The average
count values of the adjacent channel are calculated for each knot.
A fitting curve can be constructed by summing the basis spline
functions corresponding to the initial knots with the average
counts values as shown in Fig.1(b). Next, new channels are
selected in the middle of each interval as knots and the fitting
is repeated until the interval is equal to 1. With the increase of
initial knots, all the fitting curves form a finite set and the
best fitting curve can be obtained from this set according to the
criterion of Reinsch\cite{4}. The noisy spectroscopic data as
shown in Fig.2(a) is synthesized by gaussian functions with
superimposed random noise varied as a function of the square root
of the counts per channel. The fitting curve derived from the
quadric many-knot spline method is shown as Fig.2(b) corresponding
with the result using cubic B-spline shown as Fig.2(c). As can be
seen from this figure, the curve obtained by the many-knot spline
method has a similar shape as the curve obtained by B-spline
method, its validity will be described elsewhere. However, the
time consumed by both methods are entirely different as shown in
Fig.2(d). The time consumed by the many-knot spline method is
linear and changes little with the increasing of initial knots,
compared to exponentially increasing time spent using the B-spline
method in the same environment. This difference is more obvious
when the amount of spectroscopic data is large or when there are
too many initial knots.

\begin{figure}[h]
  \centering
  \subfigure[]{
    \label{fig:subfig:a}
    \includegraphics[height=6cm,width=7cm]{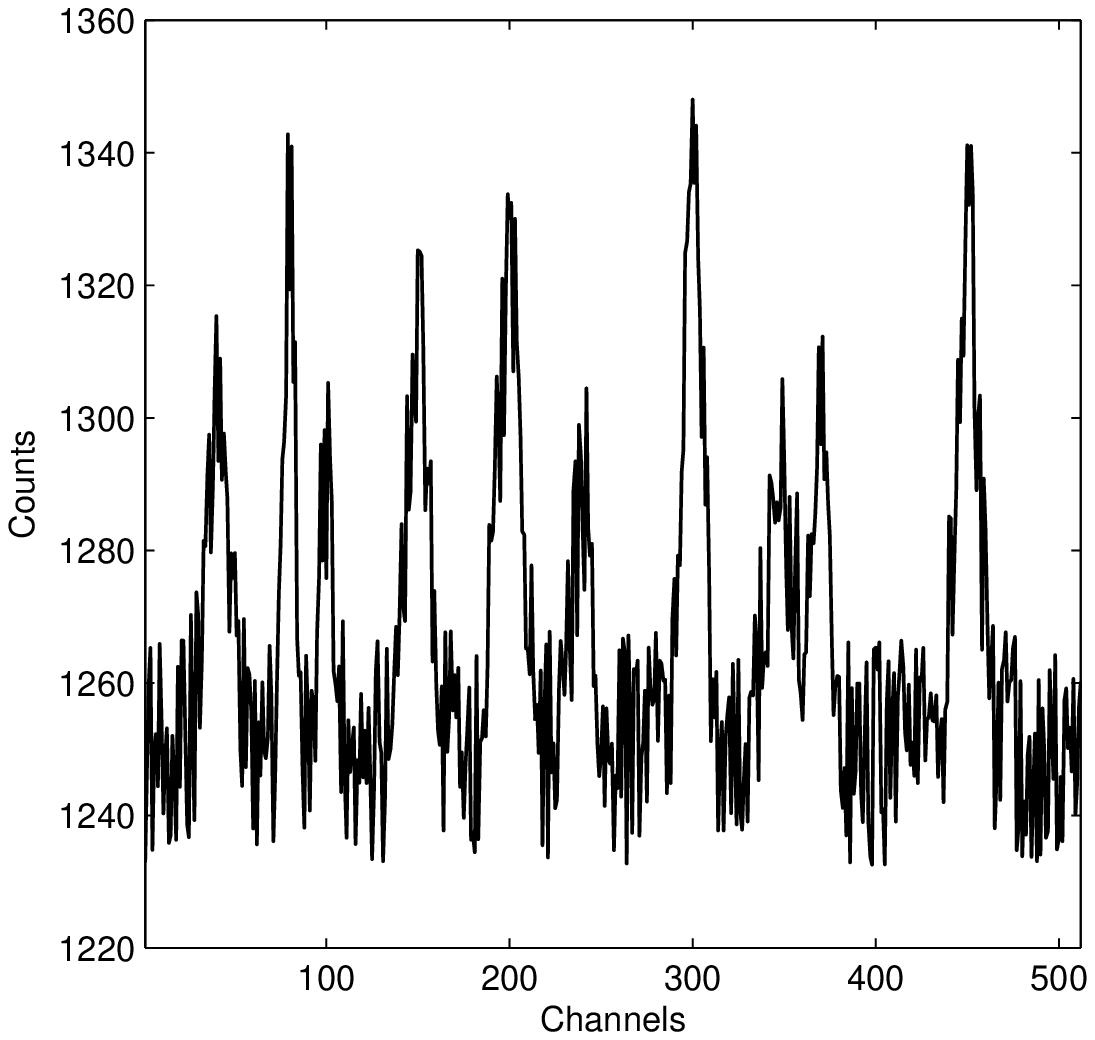}}
    \hspace{.2in}
  \subfigure[]{
    \label{fig:subfig:b}
    \includegraphics[height=6cm,width=7cm]{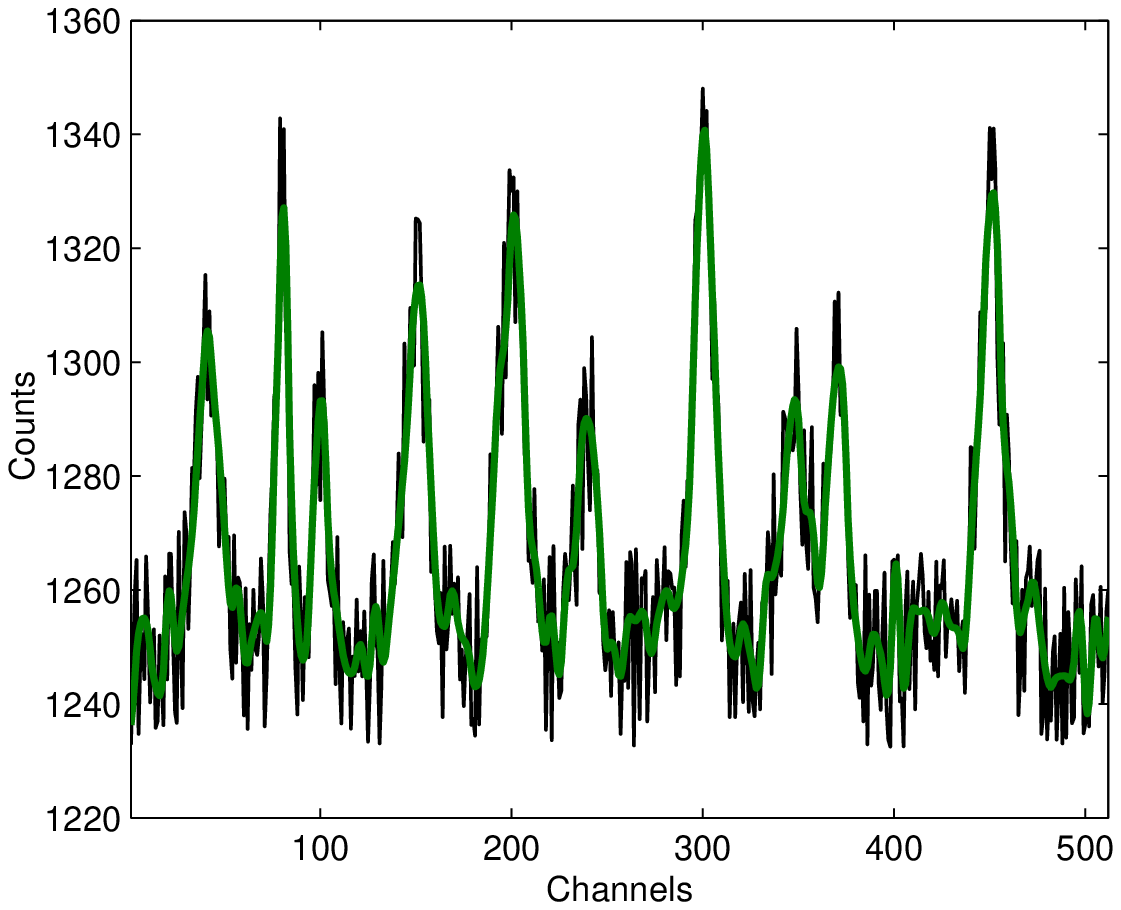}}
      \subfigure[]{
    \label{fig:subfig:a}
    \includegraphics[height=6cm,width=7cm]{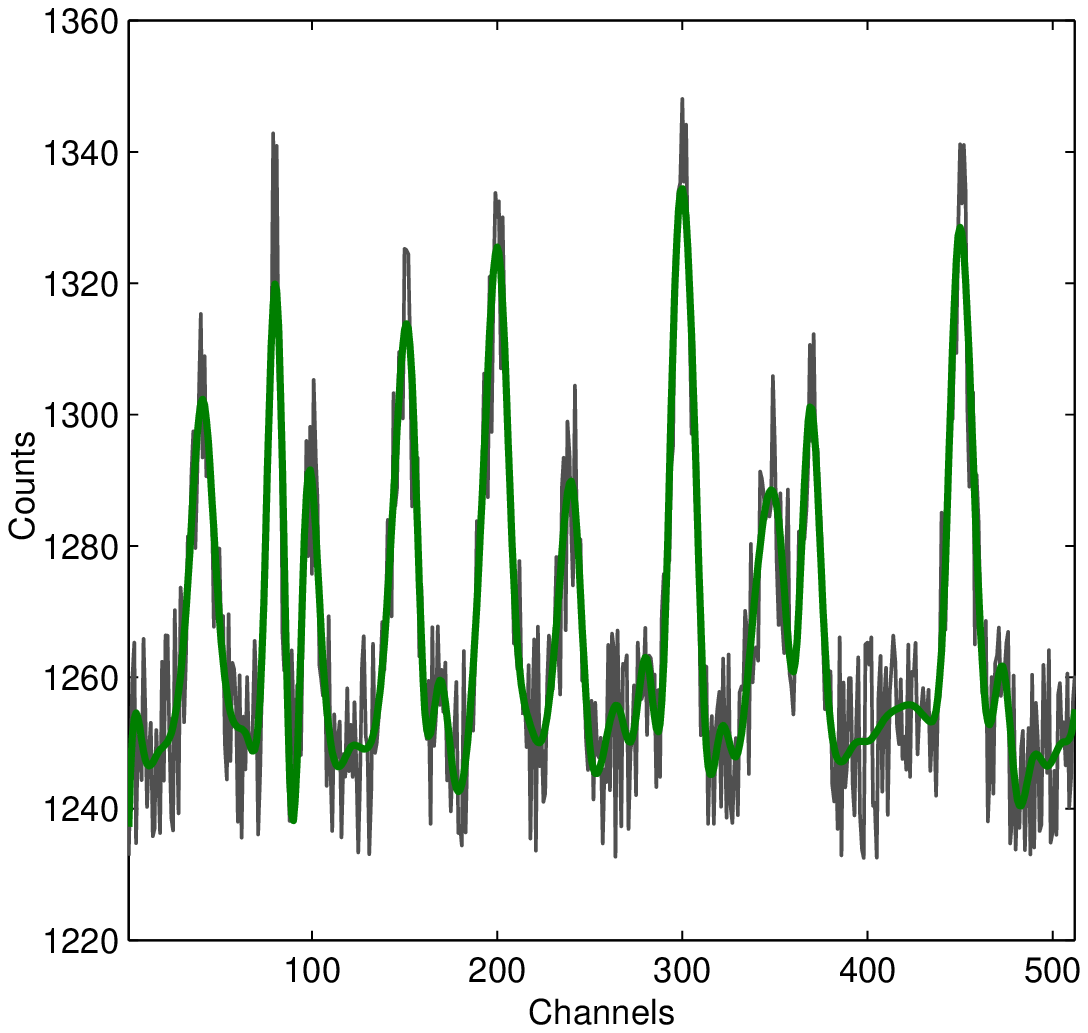}}
    \hspace{.2in}
  \subfigure[]{
    \label{fig:subfig:b}
    \includegraphics[height=6cm,width=7cm]{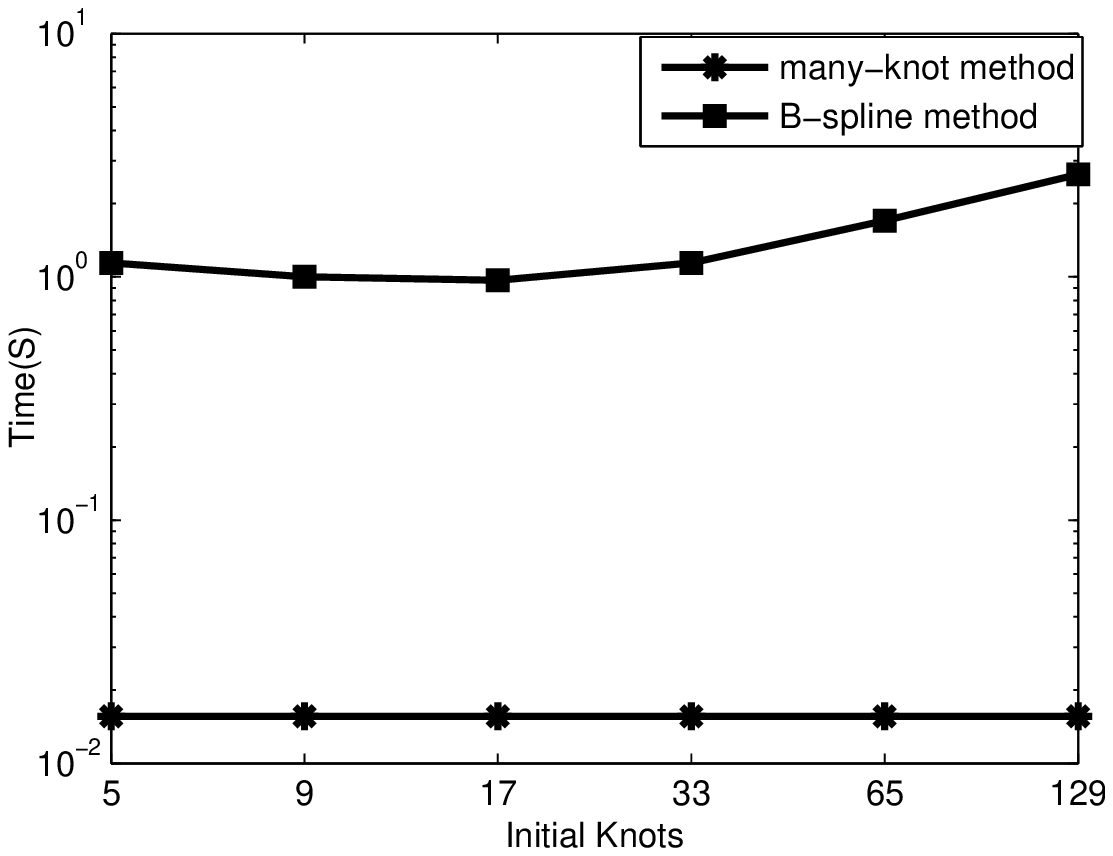}}
  \caption{(a) Noisy spectroscopic data;
(b) Smoothing noisy spectroscopic data with quadric many-knot
spline method.(c) Smoothing noisy spectroscopic data with cubic
B-spline method.(d) The time consumed by both methods with
different number of channels selected as initial knots.}
\end{figure}


\end{document}